\documentclass[12pt]{article}
\usepackage[dvipsnames]{xcolor}
\usepackage{tikz}
\usetikzlibrary{shapes}
\usetikzlibrary{tikzmark}
\tikzset{mycircled/.style={circle,draw,inner sep=0.1em,line width=0.04em}}
\usepackage{amssymb,amsmath,bm}
\usepackage{mathrsfs}
\usepackage{hyperref}
\usepackage[nomessages]{fp}
\usepackage{xparse}
\usepackage{calculator} 
\usetikzlibrary{calc} 
\hypersetup{
    colorlinks,
    citecolor=black,
    filecolor=black,
    linkcolor=black,
    urlcolor=black
}

\usepackage{color}

\newcommand{\signore}[1]{}
\usepackage{authblk}
\usepackage{scalefnt}

\usepackage{constants}
\usepackage{enumerate}
\def\p{48}
\def\minutes{15}
\def\exrate{5}
\FPeval{\minabc}{round(5*\p-12,0)}
\FPeval{\abc}{round(\p+4,0)}
\FPeval{\minsminabc}{round(\minabc*\exrate/60,1)}
\FPeval{\minsabc}{round(\abc*\exrate/60,3)}

\newcommand{\displaystyled}{
	\def\doffset{-3.2}
	\def\dscale{0.8}
	\def\ndscale{-0.8}
	\FPeval{\ndscaleovtwo}{0.5*ndscale}
	\FPeval{\dscaleovtwo}{0.5*dscale}
	\FPeval{\domright}{doffset+2*dscale}
	\FPeval{\dommid}{doffset+dscale}
	\FPeval{\legup}{dscaleovtwo*(1-0.44)}
	\FPeval{\legdown}{ndscaleovtwo*(1-0.3)}
}
\newcommand{\displaystyletwo}{
	\def\doffset{-3.2}
	\def\dscale{0.8}
	\def\ndscale{-0.8}
	\FPeval{\ndscaleovtwo}{0.5*ndscale}
	\FPeval{\dscaleovtwo}{0.5*dscale}
	\FPeval{\domright}{doffset+2*dscale}
	\FPeval{\dommid}{doffset+dscale}
	\FPeval{\legup}{dscaleovtwo*(1-0.34)}
	\FPeval{\legdown}{ndscaleovtwo*(1-0.35)}
}

\newcommand{\displaystylel}{
\def\doffset{2.6}
\def\dscale{1.1}
\def\ndscale{-1.1}
\FPeval{\ndscaleovtwo}{0.5*ndscale}
\FPeval{\dscaleovtwo}{0.5*dscale}
\FPeval{\domright}{doffset+2*dscale}
\FPeval{\dommid}{doffset+dscale}
\FPeval{\legup}{dscaleovtwo*(1-0.34)}
\FPeval{\legdown}{ndscaleovtwo*(1-0.3)} 
}

\newcommand{\inlinedom}{
\def\doffset{-0.5}
\def\dscale{0.5}
\def\ndscale{-0.5}
\FPeval{\ndscaleovtwo}{0.5*ndscale}
\FPeval{\dscaleovtwo}{0.5*dscale}
\FPeval{\domright}{doffset+2*dscale}
\FPeval{\dommid}{doffset+dscale}
\FPeval{\legup}{dscaleovtwo*(1-0.44)}
\FPeval{\legdown}{ndscaleovtwo*(1-0.3)}
}

\tikzset{
	dot hidden/.style={},
	line hidden/.style={},
	dot colour/.style={dot hidden/.append style={color=#1}},
	dot colour/.default=black,
	line colour/.style={line hidden/.append style={color=#1}},
	line colour/.default=black
}

\newcommand{\eqhinge}{2}

\newcommand{\eqhingpltwo}{4}
\newcommand{\eqhingplthree}{5}

\NewDocumentCommand{\domino}{mm}{
	\begin{tikzpicture}[x=2em,y=2em,radius=0.1]
	\draw[rounded corners=0.5,line hidden] (0,0) rectangle (2,1);
	
	\draw[line hidden] (1,0) -- (1,1);
	\ifodd#1
	\fill[dot hidden] (1.5,0.5) circle;
	\fi
	\ifnum#1>1
	\fill[dot hidden] (1.2,0.2) circle;
	\fill[dot hidden] (1.8,0.8) circle; 
	\ifnum#1>3
	\fill[dot hidden] (1.2,0.8) circle;
	\fill[dot hidden] (1.8,0.2) circle;   
	\fi
	\ifnum#1>5
	\fill[dot hidden] (1.5,0.2) circle;
	\fill[dot hidden] (1.5,0.8) circle;  
	\fi
	\fi
	\ifodd#2
	\fill[dot hidden] (0.5,0.5) circle;
	\fi
	\ifnum#2>1
	\fill[dot hidden] (0.2,0.2) circle;
	\fill[dot hidden] (0.8,0.8) circle;
	\ifnum#2>3
	\fill[dot hidden] (0.2,0.8) circle;
	\fill[dot hidden] (0.8,0.2) circle;   
	\fi
	\ifnum#2>5
	\fill[dot hidden] (0.5,0.2) circle;
	\fill[dot hidden] (0.5,0.8) circle;  
	\fi
	\fi
	\end{tikzpicture}
}

\begin{document}
\newcommand{\beq}{\begin{eqnarray}}
\newcommand{\eeq}{\end{eqnarray}}
\newcommand{\beas}{\begin{eqnarray*}}
\newcommand{\enas}{\end{eqnarray*}}
\newcommand{\bea}{\begin{eqnarray}}
\newcommand{\ena}{\end{eqnarray}}
\newcommand{\bms}{\begin{multline*}}
\newcommand{\ems}{\end{multline*}}
\newcommand{\qmq}[1]{\quad \mbox{#1} \quad}
\newcommand{\qm}[1]{\quad \mbox{#1}}
\newcommand{\nn}{\nonumber}
\newcommand{\bbox}{\hfill $\Box$}
\newcommand{\ignore}[1]{}
\newcommand{\bs}[1]{\boldsymbol{#1}}

\newcommand{\tr}{\mbox{tr}}
\newtheorem{theorem}{Theorem}[section]
\newtheorem{corollary}{Corollary}[section]
\newtheorem{conjecture}{Conjecture}[section]
\newtheorem{proposition}{Proposition}[section]
\newtheorem{remark}{Remark}[section]
\newtheorem{lemma}{Lemma}[section]
\newtheorem{definition}{Definition}[section]
\newtheorem{condition}{Condition}[section]
\newtheorem{example}{Example}[section]
\newcommand{\pf}{\noindent {\bf Proof:} }
\def\blfootnote{\xdef\@thefnmark{}\@footnotetext}

\title{The Game of Poker Chips, Dominoes and Survival}


\author{Larry Goldstein\\Department of Mathematics\\University of Southern California}

\maketitle 

\begin{abstract}
The Game of Poker Chips, Dominoes and Survival fosters team building and high level cooperation in large groups, and is a tool applied in management training exercises. Each player, initially given two colored poker chips, is allowed to make exchanges with the game coordinator according to two rules, and must secure a domino before time is called in order to `survive'. Though the rules are simple, it is not evident by their form that the survival of the entire group requires that they cooperate at a high level. From the point of view of the game coordinator, the difficulty of the game for the group can be controlled not only by the time limit, but also by the initial distribution of chips, in a way we make precise by a time complexity type argument. That analysis also provides insight into good strategies for group survival, those taking the least amount of time. In addition, coordinators may 
also want to be aware of when the game is `solvable', that is, when their initial distribution of chips permits the survival of all group members if given sufficient time to make exchanges. It turns out that the game is solvable if and only if the initial distribution contains seven chips that have one of two particular color distributions. In addition to being a lively game to play in management training or classroom settings, the analysis of the game after play can make for an engaging exercise in any basic discrete mathematics course to give a basic introduction to elements of game theory, logical reasoning, number theory and the computation of algorithmic complexities.
\end{abstract}

\section{Introduction}

Team building exercises designed to further group cohesion have long been offered to businesses and corporations as a resource for their managers and employees. Offerings run the gamut between training offered by business schools, on-line tutorials, and popular texts, see e.g.\ \cite{Dy13} or \cite{Mi15} or \cite{Ne98}. The `egg drop' exercise, for instance, consists of dividing a group into two subgroups, and, given some raw materials, having each subgroup design a package in which an uncooked egg would survive an eight foot drop. Each team is allotted time to separately discuss their ideas and come to consensus, make a pitch for their final design, and at last, see how their prototype performs in comparison to their counterpart's. Another exercise is Minefield, where the group gives verbal cues to a blindfolded member in order to navigate around obstacles. 

One somewhat mathematically oriented game that was created to promote cooperation and coordination that is currently making the rounds \cite{Ba18} uses poker chips, which come in three different colors, and dominoes. At the start of the game everyone is given two poker chips. During the game, poker chips and dominoes may be exchanged at a `bank' managed by the game coordinator, and a player must be holding a domino when time is called in order to survive. Illustrating the poker chips below as colored circles, with the number of chips of each color noted at the circle's center, the two rules for making exchanges are as follows: 

\begin{enumerate}
	\item A set of three poker chips, one of each color, 
	can be exchanged for one domino and a chip of the player's own choosing; so, for example:
	
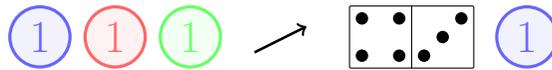
\begin{center}
		\begin{tikzpicture}
		\filldraw[color=blue!60, fill=blue!5, very thick](-1,0) circle (0.4) node {\Large{1}};
		\filldraw[color=red!60, fill=red!5, very thick](0,0) circle (0.4) node {\Large{1}};
		\filldraw[color=green!60, fill=green!5, very thick](1,0) circle (0.4) node {\Large{1}};
		\end{tikzpicture}
	 \begin{tikzpicture}[node distance=2cm]
	\pgfsetlinewidth{1pt}
	\node (A) at (0, 0) {};
	\node (B) at (1, 0.5) {};
	\draw [->] (A) edge (B);
	\end{tikzpicture}
	\domino{3}{4}
		\begin{tikzpicture}
		\filldraw[color=blue!60, fill=blue!5, very thick](-1,0) circle (0.4) node {\Large{1}};
		\end{tikzpicture}
\end{center}


\item A set of three dominoes may be exchanged for seven poker chips, all of the player's own choosing; so, for example:
	
	\begin{center}
	\domino{3}{4}
	\domino{1}{6}
	\domino{5}{2}
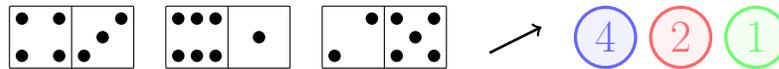
	\begin{tikzpicture}[node distance=2cm]
	\pgfsetlinewidth{1pt}
	\node (A) at (0, 0) {};
	\node (B) at (1, 0.5) {};
	\draw [->] (A) edge (B);
	\end{tikzpicture}
		\begin{tikzpicture}
		\filldraw[color=blue!60, fill=blue!5, very thick](-1,0) circle (0.4) node {\Large{4}};
		\filldraw[color=red!60, fill=red!5, very thick](0,0) circle (0.4) node {\Large{2}};
		\filldraw[color=green!60, fill=green!5, very thick](1,0) circle (0.4) node {\Large{1}};
		\end{tikzpicture}
	
	\end{center}
\end{enumerate}

Turning a group of size, let's just say \p, loose after explaining the rules and giving them two colored chips apiece and \minutes\ minutes to make the necessary exchanges, their first attempt at success may be a bit rocky, with uncoordinated bartering and bargaining, fractious subgroups, and perhaps even some chaos and mayhem. Each player, initially with only two chips, are by themselves powerless to invoke either of the two rules. Individuals who have two chips of different colors may look to make a contract with someone having a chip of the color they lack, and via Rule 1, use their pooled resources to gain a domino for the first player, while gaining a needed color for the second player. 

Due to the special properties of this game the group will soon discover that applying only such `local' exchanges will not lead to the desired state where all players hold a domino when time is called. They will, by necessity, realize that they will need to pool resources. Indeed, whereas Rule 1 can be effected by two players, even more coordination is required by Rule 2, which can only be invoked by three players all of whom are already in possession of a domino, and whose survival is therefore guaranteed. Nevertheless, Theorem \ref{thm:surv.needs.2} below shows that survival of all members of the group is impossible without invoking Rule 2; hence there must be individuals in the group who must act `selflessly' in order that all may survive. As an additional wrinkle, it may not be immediately clear how the use of Rule 2 is to anyone's advantage. Nevertheless, though experimentation and failure, the group will learn the value of Rule 2, and therefore, of cooperation.

As a practical matter, the coordinators in charge of the game need to assure that the initial distribution of colored chips to the group is sufficiently rich so that a solution exists, that is, so that the survival of all group members is possible. Clearly, some color distributions will be insufficient and the game will be futile, for instance, if all the chips distributed are of the same color. However, for five players, say, it is not immediately clear if the game has a solution when the initial distribution is 7 chips of one color, 2 of another and 1 chip of the remaining color. 

In answer to this question, Theorem \ref{thm:minsuf} shows that there exist two special sets of chips, each only of size seven, such that {\em no matter how large the size of the group}, their survival can occur if, and only if, at least one of these sets is contained in the initial distribution.  For emphasis, if only one of these special sets of size seven is distributed in a group of 1,000 players, every individual in the entire group, properly organized, and given sufficient time to make the necessary exchanges, could survive.

The reason that just one of these special sets of chips can lead to the survival of all is that it allows a `domino creation' mechanism to run, which we call the $D^3$ machine. Its inner workings are revealed by Theorem \ref{thm:D3}. Whenever survival is possible, the $D^3$ machine leads naturally to a potentially viable survival strategy. However, as time to make exchanges is limited, exclusively running the $D^3$ may not be the most efficient, that is, the quickest, way. To explore this feature of the game, in Section \ref{sec:complexity} we consider the number of exchanges required for a group to survive. We find that, depending on the initial color distribution of chips, faster strategies might be applied. Looking at the game from the point of view of the coordinator, its level of difficulty can be set from easy to impossible by controlling the colors of the chips the group initially receives.

\section{Set up and Notation} \label{sec:setup}
Given an initial distribution of chips, the game stops immediately unless the initial distribution of chips contains a `full set' of chips
\begin{tikzpicture}[scale=0.8, baseline=-1.5mm]
\filldraw[color=blue!60, fill=blue!5, very thick](-1,0) circle (0.3) node {1};
\filldraw[color=red!60, fill=red!5, very thick](-0.2,0) circle (0.3) node {1};
\filldraw[color=green!60, fill=green!5, very thick](0.6,0) circle (0.3) node {1};
\end{tikzpicture},
that is, unless the distribution contains at least one chip of each of the three colors. If it does, then the only possible move is to exchange such a triple using Rule 1. The color requested in the exchange should clearly be the least frequent one in the collection that remains, as that choice gives the most opportunities for future exchanges. The same is true for the colors returned in Rule 2, thinking of the chips returned to the players given one at a time. 

To avoid cumbersome and unnecessary bookkeeping, we can consider the chips returned in exchanges as `jokers' that can play the role of any color. The use of the joker is particularly appropriate when Rule 2 becomes frequently invoked, a sign that the group has begun coordinating. 
Including the joker chip, denoted by a circle that yet has no color, the rules of exchange now become:
\begin{enumerate}
	\item A set of three poker chips, one of each color, can be exchanged for one domino and a joker: 
	
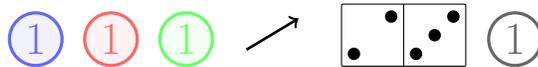
\begin{center}
		\begin{tikzpicture}
		\filldraw[color=blue!60, fill=blue!5, very thick](-1,0) circle (0.35) node {\Large{1}};
		\filldraw[color=red!60, fill=red!5, very thick](0,0) circle (0.35) node {\Large{1}};
		\filldraw[color=green!60, fill=green!5, very thick](1,0) circle (0.35) node {\Large{1}};
		\end{tikzpicture}
		\begin{tikzpicture}[node distance=2cm]
	\pgfsetlinewidth{1pt}
	\node (A) at (0, 0) {};
	\node (B) at (1, 0.6) {}; 
	\draw [->] (A) edge (B);
	\end{tikzpicture}	
	\domino{3}{2}
		\begin{tikzpicture}
		\draw[color=black!60, very thick](-1,0) circle (0.35) node {\Large{1}};
		\end{tikzpicture}
	\end{center}
	
	\item A set of three dominoes may be exchanged for seven jokers. 
	
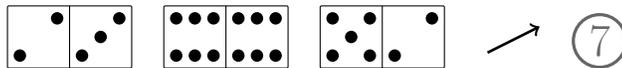
\begin{center}
\domino{3}{2}
\domino{6}{6}
\domino{2}{5}
\begin{tikzpicture}[node distance=2cm]
\pgfsetlinewidth{1pt}
\node (A) at (0, 0) {};
\node (B) at (1, 0.5) {};
\draw [->] (A) edge (B);
\end{tikzpicture}
	\begin{tikzpicture}
	\draw[color=black!60, very thick](-1,0) circle (0.35) node {\Large{7}};
	\end{tikzpicture}
\end{center}
\end{enumerate}

\displaystyled
Generally, a set of game pieces having $a,b$ and $c$ chips of the three different colors, $x$ jokers
and $d$ dominoes can be compactly visualized as
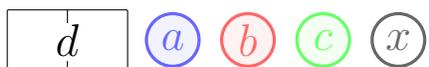
\begin{center}
	\begin{tikzpicture}
	\filldraw[color=blue!60, fill=blue!5, very thick] (-1,0) circle (0.35) node {\Large{$a$}};
	\filldraw[color=red!60, fill=red!5, very thick](0,0) circle (0.35) node {\Large{$b$}};
	\filldraw[color=green!60, fill=green!5, very thick](1,0) circle (0.35) node {\Large{$c$}};
	\draw[color=black!60, very thick](2,0) circle (0.35) node {\Large{$x$}};
	\draw[rounded corners=0.5, line hidden] (\doffset,\ndscaleovtwo) rectangle (\domright,\dscaleovtwo);
	\node at (\dommid,0) 
	{\Large{$d$}};
	\draw (\dommid,\dscaleovtwo) -- (\dommid,\legup);
	\draw (\dommid,\ndscaleovtwo) -- (\dommid,\legdown);
	\end{tikzpicture}
\end{center}
and we may denote collection of games pieces used in some exchanges in this abbreviated form. In addition, we will extend the use of the notation $\nearrow$ as used above when defining the two game rules, and will write $E \nearrow F$, read as $E$ yields $F$, if $F$ is obtained from $E$ by making successive exchanges according to either Rule 1 or Rule 2
. We will also say $E$ achieves $F$ if $E$ yields a set of which $F$ is  a subset. 

We let $I$ denote the initial set of game chips distributed to the group. We may consider arbitrary, but finite, initial distributions, and in particular those not corresponding to the game of $p$ players, where each receives two colored chips. Naturally though, that special case is of our primary interest.

As jokers can play the role of chips of any color, a consequence of Rule 1 is that all the collections 
\inlinedom
\hspace{-1.3cm}
\begin{tikzpicture}[scale=0.8, baseline=-1.5mm]
\draw[color=black!60, thick](-1,0) circle (0.3) node {3};
\end{tikzpicture}
,
\hspace{0.02cm}
\begin{tikzpicture}[scale=0.8, baseline=-1.5mm]
	\filldraw[color=blue!60, fill=blue!5, very thick](-1,0) circle (0.3) node {1};
\draw[color=black!60, thick](-0.2,0) circle (0.3) node {2};
\end{tikzpicture}
,
\hspace{0.02cm}
\begin{tikzpicture}[scale=0.8, baseline=-1.5mm]
\filldraw[color=red!60, fill=red!5, very thick](-1,0) circle (0.3) node {1};
\draw[color=black!60, thick](-0.2,0) circle (0.3) node {2};
\end{tikzpicture}
,
\hspace{0.02cm}
\begin{tikzpicture}[scale=0.8, baseline=-1.5mm]
\filldraw[color=green!60, fill=green!5, very thick](-1,0) circle (0.3) node {1};
\draw[color=black!60, thick](-0.2,0) circle (0.3) node {2};
\end{tikzpicture}
,
\hspace{0.02cm}
\begin{tikzpicture}[scale=0.8, baseline=-1.5mm]
\filldraw[color=blue!60, fill=blue!5, very thick](-1,0) circle (0.3) node {1};
\filldraw[color=red!60, fill=red!5, very thick](-0.2,0) circle (0.3) node {1};
\draw[color=black!60, thick](0.6,0) circle (0.3) node {1};
\end{tikzpicture}
,
\hspace{0.02cm}
\begin{tikzpicture}[scale=0.8, baseline=-1.5mm]
\filldraw[color=blue!60, fill=blue!5, very thick](-1,0) circle (0.3) node {1};
\filldraw[color=green!60, fill=green!5, very thick](-0.2,0) circle (0.3) node {1};
\draw[color=black!60, thick](0.6,0) circle (0.3) node {1};
\end{tikzpicture}
and
\hspace{0.02cm}
\begin{tikzpicture}[scale=0.8, baseline=-1.5mm]
\filldraw[color=red!60, fill=red!5, very thick](-1,0) circle (0.3) node {1};
\filldraw[color=green!60, fill=green!5, very thick](-0.2,0) circle (0.3) node {1};
\draw[color=black!60, thick](0.6,0) circle (0.3) node {1};
\end{tikzpicture}
yield 
\begin{tikzpicture}[scale=0.8, baseline=-1.5mm]
\draw[rounded corners=0.5, line hidden] (\doffset,\ndscaleovtwo) rectangle (\domright,\dscaleovtwo);
\node at (\dommid,0) 
{$1$};
\end{tikzpicture}
\begin{tikzpicture}[scale=0.8, baseline=-1.5mm]
\draw[color=black!60, thick](-1,0) circle (0.3) node {1};
\end{tikzpicture}.
We say we invoke Rule 1 when using Rule 1 as originally stated, or any of these consequences.

Due to the presence of the joker, choices may arise when invoking Rule 1. For instance, the set 
\hspace{0.02cm}
\begin{tikzpicture}[scale=0.8, baseline=-1.5mm]
\filldraw[color=blue!60, fill=blue!5, very thick](-1,0) circle (0.3) node {1};
\filldraw[color=red!60, fill=red!5, very thick](-0.2,0) circle (0.3) node {1};
\draw[color=black!60, thick](0.6,0) circle (0.3) node {2};
\end{tikzpicture}
can be transformed to
\begin{tikzpicture}[scale=0.8, baseline=-1.5mm]
\draw[rounded corners=0.5, line hidden] (\doffset,\ndscaleovtwo) rectangle (\domright,\dscaleovtwo);
\node at (\dommid,0) 
{$1$};
\end{tikzpicture}
\begin{tikzpicture}[scale=0.8, baseline=-1.5mm]
\filldraw[color=blue!60, fill=blue!5, very thick](-1,0) circle (0.3) node {1};
\draw[color=black!60, thick](-0.2,0) circle (0.3) node {1};
\end{tikzpicture}
by using Rule 1 on 
\begin{tikzpicture}[scale=0.8, baseline=-1.5mm]
\filldraw[color=red!60, fill=red!5, very thick](-1,0) circle (0.3) node {1};
\draw[color=black!60, thick](-0.2,0) circle (0.3) node {2};
\end{tikzpicture}
and can also be transformed to 
\begin{tikzpicture}[scale=0.8, baseline=-1.5mm]
\draw[rounded corners=0.5, line hidden] (\doffset,\ndscaleovtwo) rectangle (\domright,\dscaleovtwo);
\node at (\dommid,0) 
{$1$};
\end{tikzpicture}
\begin{tikzpicture}[scale=0.8, baseline=-1.5mm]
\draw[color=black!60, thick](-1,0) circle (0.3) node {2};
\end{tikzpicture}
by applying Rule 1 on 
\begin{tikzpicture}[scale=0.8, baseline=-1.5mm]
\filldraw[color=blue!60, fill=blue!5, very thick](-1,0) circle (0.3) node {1};
\filldraw[color=red!60, fill=red!5, very thick](-0.2,0) circle (0.3) node {1};
\draw[color=black!60, thick](0.6,0) circle (0.3) node {1};
\end{tikzpicture}.
Again, as the joker can be of any color, whatever is achievable by the first result is also achievable by the second. Generally, as Rule 1 exchanges always return a single joker, exchanges that require the least number of jokers are always preferred, leading to the imposition of the following policy that will be in force throughout. \\

\noindent {\bf Maximum Principle:} Rule 1 exchanges should always be done using sets with the fewest number of jokers. \\

It is not difficult to see that the application of the maximal principle completely determines a player's choice in any Rule 1 exchange. 
Indeed, if a Rule 1 exchange can be achieved by a given collection, than letting $r$ be the number of colors represented in the collection, we must exchange $r$ chips of the represented colors, and $3-r$ jokers, to satisfy the maximum principle.

It is simple to verify that Rule 1 and Rule 2 commute, that is, for a set in which both rules can be applied, the order in which they are is unimportant. By imposing, as we do now, that in such cases Rule 1 should be applied before Rule 2, the entire sequence of exchanges in the game become fixed once group cooperation is achieved.

When applied to a group, we use the term `survival'  as shorthand for the survival of all players. Hence, in our notation, for a game of $p$ players the desired outcome is a superset of 
\begin{tikzpicture}[scale=0.8, baseline=-1.5mm]
\draw[rounded corners=0.5, line hidden] (\doffset,\ndscaleovtwo) rectangle (\domright,\dscaleovtwo);
\node at (\dommid,0) 
{$p$};
\end{tikzpicture}.
We say a game of $p$ players is trivial if survival cannot be reached under any allowed initial distribution of chips. 

\section{Rule 1 Restricted Games}\label{sec:Rule1Restricted}

To begin to understand the fuller picture, we will first consider the game when only Rule 1 can be invoked. As the number of chips in the initial distribution $I$ is finite, and any application of Rule 1 strictly decreases the number of chips, it is clear that play using only Rule 1 must terminate, that is, reach a state where no further Rule 1 exchange is possible. We refer to a set at such a final state as a terminal set, and we say it is maximal if it achieves the maximum number of dominoes for all choices of Rule 1 exchanges on $I$.   
Given the unique path mandated by the maximum principal, and the fact that its application leads to the greatest number of exchange options, the terminal sets obtained will be maximal. 

The following result, Theorem \ref{thm:surv.needs.2}, highlights a deficit of Rule 1 with regards to the survival of the group. The intuition behind the proof starts with the observation that Rule 1, which gives a domino and a chip in return for three chips, sets the price of a domino 
at two chips. As survival costs two chips, and each individual is initially given two chips, it would seem at first glance that Rule 1 would suffice for survival, at least for some initial distributions of chips. However, for every individual to be able to exchange their two chips for a domino would require that the exchanges be efficient in the sense that all chips are converted to dominoes. But Rule 1 always returns a chip, implying that at all times at least someone in the group holds a chip. Hence, the entire collection of chips can never be entirely converted. The following result formalizes and quantifies this reasoning. 

\begin{theorem}\label{thm:surv.needs.2}
Let $n \ge 1$ denote the number of chips, of any kind or color, initially distributed in a game in which only Rule 1 may be invoked. With $d$ the number of dominoes and $y$ the number of chips at any point in the game, it holds that
\begin{align} \label{eq:preserved}
n=2d+y \qmq{with} y \ge 1, \qmq{and $y \ge 2$ when $n$ is even.}
\end{align}
No more than $p-1$ dominoes can be obtained under 
Rule 1 only play. Consequently, survival is not possible without invoking Rule 2, and all games with $p \le 3$ players are trivial. 
\end{theorem}

\noindent {\em Proof:} 
Claim \eqref{eq:preserved} holds at the start of the game with $d=0$ and $y=n$. Now say \eqref{eq:preserved} holds at some point in the game before Rule 1 is to be invoked. In that case, $y \ge 3$, as Rule 1 requires three chips. At the end of the Rule 1 exchange, these three chips are replaced by one domino and one additional chip. Hence, $y$ is replaced by $y-2$ and $d$ is replaced by $d+1$. The new value $y-2$ is therefore at least one, and the exchange leaves the sum $2d+y$ unchanged. Hence  \eqref{eq:preserved}, with $y \ge 1$, holds throughout the game. 
When $n$ is even than $2d+y$ is even, so $y$ must be even. If $y \ge 1$ and even, we must have $y \ge 2$. 

In the standard game $n=2p$ is even. As the minimal value of $y$ is 2 in this case, the maximal value of $d$ is given by the solution to $2p=2d+2$, or $d=p-1$. Hence the maximum number of dominoes that can be obtained under Rule 1 only play is $p-1$, that is, survival is not possible only invoking Rule 1, or equivalently, when never invoking Rule 2. 

To show that all games of $p \le 3$ players are trivial it suffices to show that those games can never invoke Rule 2. To invoke Rule 2 requires that $d \ge 3$, implying by \eqref{eq:preserved} that $2p = 2d+2 \ge 8$, that is, that $p \ge 4$. 
\bbox

For a game with $p$ players, Theorem \ref{thm:surv.needs.2} gives the upper bound of $p-1$ on the maximum number of dominos that can be obtained using only Rule 1. The following result show that this upper bound is achieved when $I$ consists of only jokers. We will see at the end of this section that the upper bound is also achieved for some fortuitous distributions of colored chips.

\displaystylel 
\begin{theorem} \label{thm:Xtox} For $x \ge 1$ let
\begin{align*}
	\phi(x)=\left\{
	\begin{array}{cl}
		1 & \mbox{for $x$ odd}\\2 & \mbox{for $x$ even}
	\end{array}
	\right.  
	\qmq{and}
	d(x)=\frac{x-\phi(x)}{2}.
\end{align*}	
Then, in a game where only Rule 1 may be applied, for all $x \ge 1$,
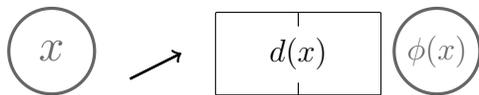
\begin{center}
	\begin{tikzpicture}
	\draw[color=black!60, very thick](-1,0) circle (0.57) node {\Large{$x$}};
	\end{tikzpicture}
	\begin{tikzpicture}[node distance=2cm]
	\pgfsetlinewidth{1pt}
	\node (A) at (0, 0) {};
	\node (B) at (1, 0.5) {};
	\draw [->] (A) edge (B);
	\end{tikzpicture}
	\begin{tikzpicture}
	\draw[rounded corners=0.5, line hidden] (\doffset,\ndscaleovtwo) rectangle (\domright,\dscaleovtwo);
	\node at (\dommid,0) 
	{$d(x)$};
	\draw (\dommid,\dscaleovtwo) -- (\dommid,\legup);
	\draw (\dommid,\ndscaleovtwo) -- (\dommid,\legdown);
	\end{tikzpicture}
	\begin{tikzpicture}
    \draw[color=black!60, very thick](-1,0) circle (0.57) node {$\phi(x)$};
	\end{tikzpicture} 
\end{center}	

This transformation is accomplished with $d(x)$ applications of Rule 1. When $2p$ jokers are distributed, $p-1$ dominoes result. 
\end{theorem}

\noindent {\em Proof:} It is clear from the definitions that $d(1)$ and $d(2)$ are both zero, and we see easily see by the structure of Rule 1 that the claim of the theorem holds for these values;  starting with two or fewer jokers, no exchanges can be made.

Consider then $x \ge 3$, when Rule 1 exchanges can be made. As 3 chips are exchanged for one domino and one chip by Rule 1, the number of chips decreases by 2, thus preserving the parity of $x$. As exchanges can continue only while the number of chips is 3 or greater, exchanges will stop when the number of chips becomes 2, or fewer. By Theorem \ref{thm:surv.needs.2}, the number of chips must be at least one. Hence, due to parity, for $x$ even, when exchanges stop there will be 2 chips, and for $x$ odd, there will be 1, that is, $\phi(x)$ chips. The number of chips expended is therefore $x-\phi(x)$, and as 2 chips are given up for each domino, the number of dominoes at stage will be half that number, that is, $d(x)$. As one domino results from each exchange, $d(x)$
is also the number of times that Rule 1 is invoked. That $d(2p)=p-1$ is
immediate from the form of $d(x)$. \bbox

To see that $p-1$ is in fact achievable in the standard game where $I$ consists of only colored chips, first uniquely write $2p=3m+r$ with $r \in \{0,1,2\}$, and consider the case where $m \ge 1$, and odd. Since $2p$ is even we must have $r=1$, and we take the initial distribution of chips to be $m+1$ of one color, and $m$ chips each of the remaining two colors. Applying Rule 1 a total of $m$ times, followed by Theorem \ref{thm:Xtox}, we obtain the upper bound $p-1$ via
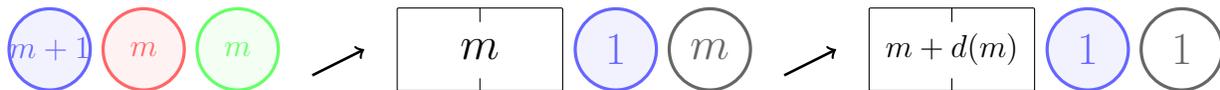
\newcommand{\circradone}{0.54}
\newcommand{\circspaceone}{1.25}
\FPeval{\cicspacetwo}{\circspaceone-1}
\FPeval{\cicspacethree}{2*\circspaceone-1}

\begin{center}
	\begin{tikzpicture}
	\filldraw[color=blue!60, fill=blue!5, very thick](-1,0) circle (\circradone) node {$m+1$};
	\filldraw[color=red!60, fill=red!5, very thick](\cicspacetwo,0) circle (\circradone) node {$m$};
	\filldraw[color=green!60, fill=green!5, very thick](\cicspacethree,0) circle (\circradone) node {$m$};
	\end{tikzpicture}
	\begin{tikzpicture}[node distance=2cm]
	\pgfsetlinewidth{1pt}
	\node (A) at (0, 0) {};
	\node (B) at (1, 0.5) {};
	\draw [->] (A) edge (B);
	\end{tikzpicture}
	\begin{tikzpicture}
	\draw[rounded corners=0.5, line hidden] (\doffset,\ndscaleovtwo) rectangle (\domright,\dscaleovtwo);
	\node at (\dommid,0) 
	{\Large{$m$}};
	\draw (\dommid,\dscaleovtwo) -- (\dommid,\legup);
	\draw (\dommid,\ndscaleovtwo) -- (\dommid,\legdown);
	\end{tikzpicture}
	\begin{tikzpicture}
	\filldraw[color=blue!60, fill=blue!5,very thick](-1,0) circle (\circradone) node {\Large{$1$}};
	\draw[color=black!60, very thick](\cicspacetwo,0) circle (\circradone) node {\Large{$m$}};
	\end{tikzpicture}
	\begin{tikzpicture}[node distance=2cm]
	\pgfsetlinewidth{1pt}
	\node (A) at (0, 0) {};
	\node (B) at (1, 0.5) {};
	\draw [->] (A) edge (B);
	\end{tikzpicture}
		\begin{tikzpicture}
	\draw[rounded corners=0.5, line hidden] (\doffset,\ndscaleovtwo) rectangle (\domright,\dscaleovtwo);
	\node at (\dommid,0) 
	{$m+d(m)$};
	\draw (\dommid,\dscaleovtwo) -- (\dommid,\legup);
	\draw (\dommid,\ndscaleovtwo) -- (\dommid,\legdown);
	\end{tikzpicture}
	\begin{tikzpicture}
	\filldraw[color=blue!60, fill=blue!5,very thick](-1,0) circle (\circradone) node {\Large{$1$}};
	\draw[color=black!60, very thick](\cicspacetwo,0) circle (\circradone) node {\Large{$1$}};
	\end{tikzpicture}
\end{center}	
using that $m$ is odd, so $\phi(m)=1$. In addition, we therefore have $d(m)=(m-1)/2$, which yields $m+d(m)=(3m-1)/2=p-1$, as desired. 

\inlinedom
We leave it to the interested reader to characterize all starting sets that achieve the maximal possible number $p-1$ of dominoes making exchanges only using Rule 1, and more generally, to determine all terminal sequences obtained from a given set, using only Rule 1. Theorem \ref{thm:Xtox}, that solves this problem for the set 
\hspace{-0.13cm}
\begin{tikzpicture}[scale=0.8, baseline=-1.5mm]
\draw[color=black!60, thick](-1,0) circle (0.3) node {$x$};
\end{tikzpicture}
gives a starting point for the exercise.

\section{The $D^3$ Machine, and Minimal Sufficient Sequences}

Section \ref{sec:Rule1Restricted} shows how applying Rule 1 alone cannot lead to success. Here we prove Theorem \ref{thm:D3}, that shows that using both rules together players can exploit a method that we will call `the $D^3$ machine'
that can, surprisingly, produce an unlimited number of domino triples.
At the machine's core lies the first set of exchanges detailed first in the theorem, which we call `joker mining', as these four exchanges produce the original collection of game pieces, with an additional joker. We show then that joker mining allows for a sequence of `domino creation' exchanges. Following this result, in Theorem \ref{thm:minsuf} we next answer the question of when the $D^3$ machine can be applied when showing that, if given sufficient time, the survival of the group can occur if and only if at least one of two special sets of chips is a subset of the initial configuration $I$. In the proof below we compress the size of dominos for typographical reasons.
\def\doffset{2.6}
\def\dscale{1.1}
\def\ndscale{-1.1}
\FPeval{\ndscaleovtwo}{0.5*ndscale}
\FPeval{\dscaleovtwo}{0.5*dscale}
\FPeval{\domright}{doffset+2*dscale}
\FPeval{\dommid}{doffset+dscale}
\FPeval{\legup}{dscaleovtwo*(1-0.34)}
\FPeval{\legdown}{ndscaleovtwo*(1-0.3)}  
\begin{theorem}[The $D^3$ machine] \label{thm:D3}
	For $r \in \{0,1\}$, joker mining 
\begin{center}
	\begin{tikzpicture}
	\draw[rounded corners=0.5, line hidden] (\doffset,\ndscaleovtwo) rectangle (\domright,\dscaleovtwo);
	\node at (\dommid,0) 
	{\Large{$3$}};
	\draw (\dommid,\dscaleovtwo) -- (\dommid,\legup);
	\draw (\dommid,\ndscaleovtwo) -- (\dommid,\legdown);
	\end{tikzpicture}
	\begin{tikzpicture}
	\draw[color=black!60, very thick](-1,0) circle (0.57) node {\Large{r}};
	\end{tikzpicture}
	\begin{tikzpicture}[node distance=2cm]
	\pgfsetlinewidth{1pt}
	\node (A) at (0, 0) {};
	\node (B) at (1, 0.5) {};
	\draw [->] (A) edge (B);
	\end{tikzpicture}
\begin{tikzpicture}
\draw[rounded corners=0.5, line hidden] (\doffset,\ndscaleovtwo) rectangle (\domright,\dscaleovtwo);
\node at (\dommid,0) 
{\Large{$3$}};
\draw (\dommid,\dscaleovtwo) -- (\dommid,\legup);
\draw (\dommid,\ndscaleovtwo) -- (\dommid,\legdown);
\end{tikzpicture}
\begin{tikzpicture}
\draw[color=black!60, very thick](-1,0) circle (0.57) node {\Large{r+1}};
\end{tikzpicture} 
\end{center}	
can be implemented in four exchanges, and domino creation 
\begin{center}
	\begin{tikzpicture} 
	\draw[rounded corners=0.5, line hidden] (\doffset,\ndscaleovtwo) rectangle (\domright,\dscaleovtwo);
	\node at (\dommid,0) 
	{\Large{$3$}};
	\draw (\dommid,\dscaleovtwo) -- (\dommid,\legup);
	\draw (\dommid,\ndscaleovtwo) -- (\dommid,\legdown);
	\end{tikzpicture}
	\begin{tikzpicture}
	\draw[color=black!60, very thick](-1,0) circle (0.57) node {\Large{1}};
	\end{tikzpicture}
	\begin{tikzpicture}[node distance=2cm]
	\pgfsetlinewidth{1pt}
	\node (A) at (0, 0) {};
	\node (B) at (1, 0.5) {};
	\draw [->] (A) edge (B);
	\end{tikzpicture}
	\begin{tikzpicture}
	\draw[rounded corners=0.5, line hidden] (\doffset,\ndscaleovtwo) rectangle (\domright,\dscaleovtwo);
	\node at (\dommid,0) 
	{\Large{$4$}};
	\draw (\dommid,\dscaleovtwo) -- (\dommid,\legup);
	\draw (\dommid,\ndscaleovtwo) -- (\dommid,\legdown);
	\end{tikzpicture}
	\begin{tikzpicture}
	\draw[color=black!60, very thick](-1,0) circle (0.57) node {\Large{1}};
	\end{tikzpicture}
\end{center}	
implemented by nine.
For all $p \ge 3$ we have 
\def\doffset{-0.5}
\def\dscale{0.5}
\def\ndscale{-0.5}
\FPeval{\ndscaleovtwo}{0.5*ndscale}
\FPeval{\dscaleovtwo}{0.5*dscale}
\FPeval{\domright}{doffset+2*dscale}
\FPeval{\dommid}{doffset+dscale}
\FPeval{\legup}{dscaleovtwo*(1-0.44)}
\FPeval{\legdown}{ndscaleovtwo*(1-0.3)}
\hspace{-1.3cm}
\begin{tikzpicture}[scale=0.8, baseline=-1.5mm]
\draw[rounded corners=0.5, line hidden] (\doffset,\ndscaleovtwo) rectangle (\domright,\dscaleovtwo);
\node at (\dommid,0) 
{$3$};
\end{tikzpicture}
$\nearrow$
\begin{tikzpicture}[scale=0.8, baseline=-1.5mm]
\draw[rounded corners=0.5, line hidden] (\doffset,\ndscaleovtwo) rectangle (\domright,\dscaleovtwo);
\node at (\dommid,0) 
{$p$};
\end{tikzpicture},
and a group of size $p$ can achieve 
\begin{tikzpicture}[scale=0.8, baseline=-1.5mm]
\draw[rounded corners=0.5, line hidden] (\doffset,\ndscaleovtwo) rectangle (\domright,\dscaleovtwo);
\node at (\dommid,0) 
{$p$};
\end{tikzpicture}
if and only if $I$ achieves
\begin{tikzpicture}[scale=0.8, baseline=-1.5mm]
\draw[rounded corners=0.5, line hidden] (\doffset,\ndscaleovtwo) rectangle (\domright,\dscaleovtwo);
\node at (\dommid,0) 
{$3$};
\end{tikzpicture}.
\end{theorem}

\noindent {\em Proof:}
Applying Rule 2 on
\def\doffset{-0.5}
\def\dscale{0.5}
\def\ndscale{-0.5}
\FPeval{\ndscaleovtwo}{0.5*ndscale}
\FPeval{\dscaleovtwo}{0.5*dscale}
\FPeval{\domright}{doffset+2*dscale}
\FPeval{\dommid}{doffset+dscale}
\FPeval{\legup}{dscaleovtwo*(1-0.44)}
\FPeval{\legdown}{ndscaleovtwo*(1-0.3)}
\hspace{-1.3cm}
\begin{tikzpicture}[scale=0.8, baseline=-1.5mm]
\draw[rounded corners=0.5, line hidden] (\doffset,\ndscaleovtwo) rectangle (\domright,\dscaleovtwo);
\node at (\dommid,0) 
{$3$};
\end{tikzpicture}
\begin{tikzpicture}[scale=0.8, baseline=-1.5mm]
\draw[color=black!60, thick](-1,0) circle (0.3) node {$r$};
\end{tikzpicture}
for the first exchange, and then Rule 1 for the remaining three exchanges, yields
\def\scale{1}
\def\aspect{1.2}
\def\domleft{-0.5}
\def\lowery{-0.15} 
\def\ulegscale{0.12} 
\def\llegscale{0.12}  
\def\temparrow{0.8}
\FPeval{\temprad}{0.5*scale}
\FPeval{\domright}{domleft+aspect*scale}
\FPeval{\uppery}{lowery+scale}
\FPeval{\midx}{domleft+0.5*aspect*scale}
\FPeval{\midy}{lowery+0.5*scale}
\FPeval{\legup}{uppery-scale*ulegscale}
\FPeval{\legdown}{lowery+scale*llegscale}
\begin{center}
\begin{tikzpicture}
\draw[rounded corners=0.5, line hidden] 
(\domleft,\lowery) rectangle (\domright,\uppery);
\node at (\midx,\midy) 
{\Large{3}};
\draw (\midx,\uppery) -- (\midx,\legup);
\draw (\midx,\lowery) -- (\midx,\legdown);
\end{tikzpicture}
\begin{tikzpicture}
\draw[color=black!60, very thick](-1,0) circle (\temprad) node {\Large{$r$}};
\end{tikzpicture}
\begin{tikzpicture}[node distance=2cm]
\pgfsetlinewidth{1pt}
\node (A) at (0, 0) {};
\node (B) at (\temparrow, 0.5) {};
\draw [->] (A) edge (B);
\end{tikzpicture}
\begin{tikzpicture}
\draw[color=black!60, very thick](-1,0) circle (\temprad) node {$r+7$};
\end{tikzpicture}
\begin{tikzpicture}[node distance=2cm]
\pgfsetlinewidth{1pt}
\node (A) at (0, 0) {};
\node (B) at (\temparrow, 0.5) {};
\draw [->] (A) edge (B);
\end{tikzpicture}
\begin{tikzpicture}
\draw[rounded corners=0.5, line hidden] 
(\domleft,\lowery) rectangle (\domright,\uppery);
\node at (\midx,\midy) 
{\Large{1}};
\draw (\midx,\uppery) -- (\midx,\legup);
\draw (\midx,\lowery) -- (\midx,\legdown);
\end{tikzpicture}
\begin{tikzpicture}
\draw[color=black!60, very thick](-1,0) circle (\temprad) node {$r+5$};
\end{tikzpicture}
\begin{tikzpicture}[node distance=2cm]
\pgfsetlinewidth{1pt}
\node (A) at (0, 0) {};
\node (B) at (\temparrow, 0.5) {};
\draw [->] (A) edge (B);
\end{tikzpicture}
\begin{tikzpicture}
\draw[rounded corners=0.5, line hidden] 
(\domleft,\lowery) rectangle (\domright,\uppery);
\node at (\midx,\midy) 
{\Large{2}};
\draw (\midx,\uppery) -- (\midx,\legup);
\draw (\midx,\lowery) -- (\midx,\legdown);
\end{tikzpicture}
\begin{tikzpicture}
\draw[color=black!60, very thick](-1,0) circle (\temprad) node {$r+3$};
\end{tikzpicture}
\begin{tikzpicture}[node distance=2cm]
\pgfsetlinewidth{1pt}
\node (A) at (0, 0) {};
\node (B) at (\temparrow, 0.5) {};
\draw [->] (A) edge (B);
\end{tikzpicture}
\begin{tikzpicture}
\draw[rounded corners=0.5, line hidden] 
(\domleft,\lowery) rectangle (\domright,\uppery);
\node at (\midx,\midy) 
{\Large{3}};
\draw (\midx,\uppery) -- (\midx,\legup);
\draw (\midx,\lowery) -- (\midx,\legdown);
\end{tikzpicture}
\begin{tikzpicture}
\draw[color=black!60, very thick](-1,0) circle (\temprad) node {$r+1$};
\end{tikzpicture}
\end{center}
showing the first claim.

Applying joker mining twice, at a cost of 8 exchanges, we obtain
\begin{center}
	\begin{tikzpicture}
	\draw[rounded corners=0.5, line hidden] 
	(\domleft,\lowery) rectangle (\domright,\uppery);
	\node at (\midx,\midy) 
	{\Large{3}};
	\draw (\midx,\uppery) -- (\midx,\legup);
	\draw (\midx,\lowery) -- (\midx,\legdown);
	\end{tikzpicture}
	\begin{tikzpicture}
	\draw[color=black!60, very thick](-1,0) circle (\temprad) node {\Large{$r$}};
	\end{tikzpicture}
	\begin{tikzpicture}[node distance=2cm]
	\pgfsetlinewidth{1pt}
	\node (A) at (0, 0) {};
	\node (B) at (\temparrow, 0.5) {};
	\draw [->] (A) edge (B);
	\end{tikzpicture}
	\begin{tikzpicture}
	\draw[rounded corners=0.5, line hidden] 
	(\domleft,\lowery) rectangle (\domright,\uppery);
	\node at (\midx,\midy) 
	{\Large{3}};
	\draw (\midx,\uppery) -- (\midx,\legup);
	\draw (\midx,\lowery) -- (\midx,\legdown);
	\end{tikzpicture}
	\begin{tikzpicture}
	\draw[color=black!60, very thick](-1,0) circle (\temprad) node {$r+1$};
	\end{tikzpicture}
	\begin{tikzpicture}[node distance=2cm]
	\pgfsetlinewidth{1pt}
	\node (A) at (0, 0) {};
	\node (B) at (\temparrow, 0.5) {};
	\draw [->] (A) edge (B);
	\end{tikzpicture}
	\begin{tikzpicture}
	\draw[rounded corners=0.5, line hidden] 
	(\domleft,\lowery) rectangle (\domright,\uppery);
	\node at (\midx,\midy) 
	{\Large{3}};
	\draw (\midx,\uppery) -- (\midx,\legup);
	\draw (\midx,\lowery) -- (\midx,\legdown);
	\end{tikzpicture}
	\begin{tikzpicture}
	\draw[color=black!60, very thick](-1,0) circle (\temprad) node {$r+2$};
	\end{tikzpicture}
\end{center}

When $r=1$ applying Rule 1 once to the result of this sequence of yields, for a total of 9 exchanges, the sequence
\begin{center}
	\begin{tikzpicture}
	\draw[rounded corners=0.5, line hidden] 
	(\domleft,\lowery) rectangle (\domright,\uppery);
	\node at (\midx,\midy) 
	{\Large{3}};
	\draw (\midx,\uppery) -- (\midx,\legup);
	\draw (\midx,\lowery) -- (\midx,\legdown);
	\end{tikzpicture}
	\begin{tikzpicture}
	\draw[color=black!60, very thick](-1,0) circle (\temprad) node {\Large{$1$}};
	\end{tikzpicture}
	\begin{tikzpicture}[node distance=2cm]
	\pgfsetlinewidth{1pt}
	\node (A) at (0, 0) {};
	\node (B) at (\temparrow, 0.5) {};
	\draw [->] (A) edge (B);
	\end{tikzpicture}
	\begin{tikzpicture}
	\draw[rounded corners=0.5, line hidden] 
	(\domleft,\lowery) rectangle (\domright,\uppery);
	\node at (\midx,\midy) 
	{\Large{3}};
	\draw (\midx,\uppery) -- (\midx,\legup);
	\draw (\midx,\lowery) -- (\midx,\legdown);
	\end{tikzpicture}
	\begin{tikzpicture}
	\draw[color=black!60, very thick](-1,0) circle (\temprad) node {\Large{$3$}};
	\end{tikzpicture}
	\begin{tikzpicture}[node distance=2cm]
	\pgfsetlinewidth{1pt}
	\node (A) at (0, 0) {};
	\node (B) at (\temparrow, 0.5) {};
	\draw [->] (A) edge (B);
	\end{tikzpicture}
	\begin{tikzpicture}
	\draw[rounded corners=0.5, line hidden] 
	(\domleft,\lowery) rectangle (\domright,\uppery);
	\node at (\midx,\midy) 
	{\Large{4}};
	\draw (\midx,\uppery) -- (\midx,\legup);
	\draw (\midx,\lowery) -- (\midx,\legdown);
	\end{tikzpicture}
	\begin{tikzpicture}
	\draw[color=black!60, very thick](-1,0) circle (\temprad) node {\Large{$1$}};
	\end{tikzpicture}
\end{center}
as claimed. When $r=0$, 
\ignore{and applying the transformation \eqref{eq:Xmining} to the result of \eqref{eq:D3Xr2} one additional time when $r=0$,  costing 4 additional exchanges for a total of 12, followed by Rule 1, we have
$$
D^3X^2 \rightarrow D^3X^3 \rightarrow  D^4X,
$$
and hence $D^3 \rightarrow D^4X$ at a cost of 13 total exchanges. }
applying joker mining three times followed by a Rule 1 exchange produces 
\begin{center}
	\begin{tikzpicture}
	\draw[rounded corners=0.5, line hidden] 
	(\domleft,\lowery) rectangle (\domright,\uppery);
	\node at (\midx,\midy) 
	{\Large{3}};
	\draw (\midx,\uppery) -- (\midx,\legup);
	\draw (\midx,\lowery) -- (\midx,\legdown);
	\end{tikzpicture}
	\begin{tikzpicture}[node distance=2cm]
	\pgfsetlinewidth{1pt}
	\node (A) at (0, 0) {};
	\node (B) at (\temparrow, 0.5) {};
	\draw [->] (A) edge (B);
	\end{tikzpicture}
	\begin{tikzpicture}
	\draw[rounded corners=0.5, line hidden] 
	(\domleft,\lowery) rectangle (\domright,\uppery);
	\node at (\midx,\midy) 
	{\Large{3}};
	\draw (\midx,\uppery) -- (\midx,\legup);
	\draw (\midx,\lowery) -- (\midx,\legdown);
	\end{tikzpicture}
	\begin{tikzpicture}
	\draw[color=black!60, very thick](-1,0) circle (\temprad) node {\Large{$3$}};
	\end{tikzpicture}
	\begin{tikzpicture}[node distance=2cm]
	\pgfsetlinewidth{1pt}
	\node (A) at (0, 0) {};
	\node (B) at (\temparrow, 0.5) {};
	\draw [->] (A) edge (B);
	\end{tikzpicture}
	\begin{tikzpicture}
	\draw[rounded corners=0.5, line hidden] 
	(\domleft,\lowery) rectangle (\domright,\uppery);
	\node at (\midx,\midy) 
	{\Large{4}};
	\draw (\midx,\uppery) -- (\midx,\legup);
	\draw (\midx,\lowery) -- (\midx,\legdown);
	\end{tikzpicture}
	\begin{tikzpicture}
	\draw[color=black!60, very thick](-1,0) circle (\temprad) node {\Large{$1$}};
	\end{tikzpicture}
\end{center} 
Clearly any number of $p \ge 3$ dominoes can 
can be produced by repeating this sequence of exchanges as necessary.

For the final claim, we have already shown that if $p \ge 4$ and $I$ achieves 
\begin{tikzpicture}[scale=0.8, baseline=-1.5mm]
\draw[rounded corners=0.5, line hidden] (\doffset,\ndscaleovtwo) rectangle (\domright,\dscaleovtwo);
\node at (\dommid,0) 
{$3$};
\end{tikzpicture}
then 
\begin{tikzpicture}[scale=0.8, baseline=-1.5mm]
\draw[rounded corners=0.5, line hidden] (\doffset,\ndscaleovtwo) rectangle (\domright,\dscaleovtwo);
\node at (\dommid,0) 
{$p$};
\end{tikzpicture}
can be achieved. Conversely, if $I$ does not achieve 
\begin{tikzpicture}[scale=0.8, baseline=-1.5mm]
\draw[rounded corners=0.5, line hidden] (\doffset,\ndscaleovtwo) rectangle (\domright,\dscaleovtwo);
\node at (\dommid,0) 
{$3$};
\end{tikzpicture}
then Rule 2 can never be invoked, and Theorem \ref{thm:surv.needs.2} shows that survival is not possible. 
\bbox\\

Joker mining may at first appear to produce a game piece from thin air, taking as input 
\begin{tikzpicture}[scale=0.8, baseline=-1.5mm]
\draw[rounded corners=0.5, line hidden] (\doffset,\ndscaleovtwo) rectangle (\domright,\dscaleovtwo);
\node at (\dommid,0) 
{$3$};
\end{tikzpicture}
and yielding
\begin{tikzpicture}[scale=0.8, baseline=-1.5mm]
\draw[rounded corners=0.5, line hidden] (\doffset,\ndscaleovtwo) rectangle (\domright,\dscaleovtwo);
\node at (\dommid,0) 
{$3$};
\end{tikzpicture}
\begin{tikzpicture}[scale=0.8, baseline=-1.5mm]
\draw[color=black!60, thick](-1,0) circle (0.3) node {1};
\end{tikzpicture}.
However, if we recall that the cost of a domino, as fixed by Rule 1, is two chips, the exchange of three dominoes from seven chips reveals the source of the benefit.

We now turn to the question of which initial sets lead to the survival of all group members, given enough time to make exchanges. To take an extreme case, clearly a group cannot survive if they are initially given chips all of the same color.  We say that the set $E$ of game pieces is sufficient for
\begin{tikzpicture}[scale=0.8, baseline=-1.5mm]
\draw[rounded corners=0.5, line hidden] (\doffset,\ndscaleovtwo) rectangle (\domright,\dscaleovtwo);
\node at (\dommid,0) 
{$3$};
\end{tikzpicture},
or just sufficient, if $E$ achieves
\begin{tikzpicture}[scale=0.8, baseline=-1.5mm]
\draw[rounded corners=0.5, line hidden] (\doffset,\ndscaleovtwo) rectangle (\domright,\dscaleovtwo);
\node at (\dommid,0) 
{$3$};
\end{tikzpicture}.
We say $E$ is minimal sufficient if $E$ is sufficient and no subset of $E$ achieves 
\begin{tikzpicture}[scale=0.8, baseline=-1.5mm]
\draw[rounded corners=0.5, line hidden] (\doffset,\ndscaleovtwo) rectangle (\domright,\dscaleovtwo);
\node at (\dommid,0) 
{$3$};
\end{tikzpicture}.

\begin{theorem}\label{thm:minsuf}
Up to a permutation of colors, there are exactly two minimal sufficient sets, 
\begin{tikzpicture}[scale=0.8, baseline=-1.5mm]
\filldraw[color=blue!60, fill=blue!5, very thick](-1,0) circle (0.3) node {3};
\filldraw[color=red!60, fill=red!5, very thick](-0.2,0) circle (0.3) node {3};
\filldraw[color=green!60, fill=green!5, very thick](0.6,0) circle (0.3) node {1};
\end{tikzpicture}
and 
\begin{tikzpicture}[scale=0.8, baseline=-1.5mm]
\filldraw[color=blue!60, fill=blue!5, very thick](-1,0) circle (0.3) node {3};
\filldraw[color=red!60, fill=red!5, very thick](-0.2,0) circle (0.3) node {2};
\filldraw[color=green!60, fill=green!5, very thick](0.6,0) circle (0.3) node {2};
\end{tikzpicture}, 
and for a group given sufficient time to make exchanges, survival can occur if and only if at least one of these sets is a subset of $I$.  Games with $p$ players are non-trivial if and only if $p \ge 4$. 
\end{theorem}

\noindent {\em Proof:} By relabeling the colors, we may assume without loss of generality that the intital distribution of chips is given by $I=$
\begin{tikzpicture}[scale=0.8, baseline=-1.5mm]
\filldraw[color=blue!60, fill=blue!5, very thick](-1,0) circle (0.3) node {$a$};
\filldraw[color=red!60, fill=red!5, very thick](-0.2,0) circle (0.3) node {$b$};
\filldraw[color=green!60, fill=green!5, very thick](0.6,0) circle (0.3) node {$c$};
\end{tikzpicture},
with $a \ge b \ge c$. No exchanges can be made unless $c \ge 1$, in which case only Rule 1 may be invoked, leading to

\hspace{1.2cm}
\begin{tikzpicture}
	\filldraw[color=blue!60, fill=blue!5, very thick](-1,0) circle (0.35) node {a};
	\filldraw[color=red!60, fill=red!5, very thick](0,0) circle (0.35) node {b};
	\filldraw[color=green!60, fill=green!5, very thick](1,0) circle (0.35) node {c};
	\end{tikzpicture}
	\begin{tikzpicture}[node distance=2cm]
	\pgfsetlinewidth{1pt}
	\node (A) at (0, 0) {};
	\node (B) at (1, 0.6) {};
	\draw [->] (A) edge (B);
	\end{tikzpicture}
	\domino{3}{6}
	\begin{tikzpicture}
	\filldraw[color=blue!60, fill=blue!5, very thick](-1,0) circle (0.35) node {a-1};
	\filldraw[color=red!60, fill=red!5, very thick](0,0) circle (0.35) node {b-1};
	\filldraw[color=green!60, fill=green!5, very thick](1,0) circle (0.35) node {c-1};
	\end{tikzpicture}
	\begin{tikzpicture}
	\draw[color=black!60, very thick](-1,0) circle (0.35) node {1};
	\node[] at (3.6,0) {(\eqhinge)};
	\end{tikzpicture}

First consider the case $c=1$. No further exchanges can result from (\eqhinge) unless $b \ge 2$, and then the only option available is to apply Rule 1 to obtain
\begin{center}
	\domino{3}{6}
	\begin{tikzpicture}
	\filldraw[color=blue!60, fill=blue!5, very thick](-1,0) circle (0.35) node {a-1};
	\filldraw[color=red!60, fill=red!5, very thick](0,0) circle (0.35) node {b-1};
	\end{tikzpicture}
	\begin{tikzpicture}
	\draw[color=black!60, very thick](-1,0) circle (0.35) node {1};
	\end{tikzpicture}
	\begin{tikzpicture}[node distance=2cm]
	\pgfsetlinewidth{1pt}
	\node (A) at (0, 0) {};
	\node (B) at (1, 0.6) {};
	\draw [->] (A) edge (B);
	\end{tikzpicture}
	\domino{3}{6}
	\domino{1}{4}
	\begin{tikzpicture}
	\filldraw[color=blue!60, fill=blue!5, very thick](-1,0) circle (0.35) node {a-2};
	\filldraw[color=red!60, fill=red!5, very thick](0,0) circle (0.35) node {b-2};
	\end{tikzpicture}
	\begin{tikzpicture}
	\draw[color=black!60, very thick](-1,0) circle (0.35) node {1};
	\end{tikzpicture}
\end{center}

Now we must impose the condition that $b \ge 3$ in order to continue, in which case, again using Rule 1, we achieve the set
\begin{center}
	\domino{2}{5}
	\domino{3}{6}
	\domino{1}{4}
	\begin{tikzpicture}
	\filldraw[color=blue!60, fill=blue!5, very thick](-1,0) circle (0.35) node {a-3};
	\filldraw[color=red!60, fill=red!5, very thick](0,0) circle (0.35) node {b-3};
	\end{tikzpicture}
	\begin{tikzpicture}
	\draw[color=black!60, very thick](-1,0) circle (0.35) node {1};
	\end{tikzpicture}
\end{center}
Hence, the minimal set that can achieve 
\begin{tikzpicture}[scale=0.8, baseline=-1.5mm]
\draw[rounded corners=0.5, line hidden] (\doffset,\ndscaleovtwo) rectangle (\domright,\dscaleovtwo);
\node at (\dommid,0) 
{$3$};
\end{tikzpicture}
when $c=1$ is the set $M=$
\begin{tikzpicture}[scale=0.8, baseline=-1.5mm]
\filldraw[color=blue!60, fill=blue!5, very thick](-1,0) circle (0.3) node {3};
\filldraw[color=red!60, fill=red!5, very thick](-0.2,0) circle (0.3) node {3};
\filldraw[color=green!60, fill=green!5, very thick](0.6,0) circle (0.3) node {1};
\end{tikzpicture}.

In the remaining case $c \ge 2$, and starting from the right hand side of (\eqhinge), the maximal principle dictates the exchange
\begin{center}
\domino{3}{6}
\begin{tikzpicture}
\filldraw[color=blue!60, fill=blue!5, very thick](-1,0) circle (0.35) node {a-1};
\filldraw[color=red!60, fill=red!5, very thick](0,0) circle (0.35) node {b-1};
\filldraw[color=green!60, fill=green!5, very thick](1,0) circle (0.35) node {c-1};
\end{tikzpicture}
\begin{tikzpicture}
\draw[color=black!60, very thick](-1,0) circle (0.35) node {1};
\end{tikzpicture}
\begin{tikzpicture}[node distance=2cm]
\pgfsetlinewidth{1pt}
\node (A) at (0, 0) {};
\node (B) at (1, 0.6) {};
\draw [->] (A) edge (B);
\end{tikzpicture}
\domino{3}{6}
\domino{1}{4}
\begin{tikzpicture}
\filldraw[color=blue!60, fill=blue!5, very thick](-1,0) circle (0.35) node {a-2};
\filldraw[color=red!60, fill=red!5, very thick](0,0) circle (0.35) node {b-2};
\filldraw[color=green!60, fill=green!5, very thick](1,0) circle (0.35) node {c-2};
\end{tikzpicture}
\begin{tikzpicture}
\draw[color=black!60, very thick](-1,0) circle (0.35) node {2};
\end{tikzpicture}
\end{center}

At this point, we must impose the condition that $a \ge 3$ for the process to continue. Initial sets $I$ with $b \ge 3$ will contain 
\begin{tikzpicture}[scale=0.8, baseline=-1.5mm]
\filldraw[color=blue!60, fill=blue!5, very thick](-1,0) circle (0.3) node {3};
\filldraw[color=red!60, fill=red!5, very thick](-0.2,0) circle (0.3) node {3};
\filldraw[color=green!60, fill=green!5, very thick](0.6,0) circle (0.3) node {1};
\end{tikzpicture}
as a proper subset, so will achieve
\begin{tikzpicture}[scale=0.8, baseline=-1.5mm]
\draw[rounded corners=0.5, line hidden] (\doffset,\ndscaleovtwo) rectangle (\domright,\dscaleovtwo);
\node at (\dommid,0) 
{$3$};
\end{tikzpicture},
but will not be minimal. Hence, we are left with the constraints $a \ge 3, b=2,c=2$, for which we may continue by making the exchange

\displaystyletwo 
\begin{center}
\begin{tikzpicture}
\draw[rounded corners=0.5, line hidden] (\doffset,\ndscaleovtwo) rectangle (\domright,\dscaleovtwo);
\node at (\dommid,0) 
{\Large{$2$}};
\draw (\dommid,\dscaleovtwo) -- (\dommid,\legup);
\draw (\dommid,\ndscaleovtwo) -- (\dommid,\legdown);
\filldraw[color=blue!60, fill=blue!5, very thick] (-1,0) circle (0.35) node {a-2};
\draw[color=black!60, very thick](0,0) circle (0.35) node {2};
\end{tikzpicture}
	\begin{tikzpicture}[node distance=2cm]
	\pgfsetlinewidth{1pt}
	\node (A) at (0, 0) {};
	\node (B) at (1, 0.6) {};
	\draw [->] (A) edge (B);
	\end{tikzpicture}
\begin{tikzpicture}
\draw[rounded corners=0.5, line hidden] (\doffset,\ndscaleovtwo) rectangle (\domright,\dscaleovtwo);
\node at (\dommid,0) 
{\Large{$3$}};
\draw (\dommid,\dscaleovtwo) -- (\dommid,\legup);
\draw (\dommid,\ndscaleovtwo) -- (\dommid,\legdown);
\filldraw[color=blue!60, fill=blue!5, very thick] (-1,0) circle (0.35) node {a-3};
\draw[color=black!60, very thick](0,0) circle (0.35) node {1};
\end{tikzpicture}
\end{center}
Hence, the minimal set that achieves 
\inlinedom
\hspace{-1.93cm}
\begin{tikzpicture}[scale=0.8, baseline=-1.5mm]
\draw[rounded corners=0.5, line hidden] (\doffset,\ndscaleovtwo) rectangle (\domright,\dscaleovtwo);
\node at (\dommid,0) 
{$3$};
\end{tikzpicture}
when $c \ge 2$ is $N=$
\begin{tikzpicture}[scale=0.8, baseline=-1.5mm]
\filldraw[color=blue!60, fill=blue!5, very thick](-1,0) circle (0.3) node {3};
\filldraw[color=red!60, fill=red!5, very thick](-0.2,0) circle (0.3) node {2};
\filldraw[color=green!60, fill=green!5, very thick](0.6,0) circle (0.3) node {2};
\end{tikzpicture}.
As the cases considered are exhaustive, $M$ and $N$ are the only minimial sufficient sets.

Whenever one of $M$ or $N$ is a subset of $I$, then 
\begin{tikzpicture}[scale=0.8, baseline=-1.5mm]
\draw[rounded corners=0.5, line hidden] (\doffset,\ndscaleovtwo) rectangle (\domright,\dscaleovtwo);
\node at (\dommid,0) 
{$3$};
\end{tikzpicture}
can be achieved and survival can occur by Theorem \ref{thm:D3}. Conversely, if neither $M$ nor $N$ is a subset of $I$, and survival occurs, then $I$ achieves
\begin{tikzpicture}[scale=0.8, baseline=-1.5mm]
\draw[rounded corners=0.5, line hidden] (\doffset,\ndscaleovtwo) rectangle (\domright,\dscaleovtwo);
\node at (\dommid,0) 
{$3$};
\end{tikzpicture}
by Theorem \ref{thm:surv.needs.2}. In that case, $I$, or some subset of $I$, is a minimal sufficient set that does not contain $M$ or $N$ as a subset, which is a contradiction.

For the final claim on non-triviality, Theorem \ref{thm:surv.needs.2} gives that games with $p \le 3$ players are trivial. The converse obtains noting that the minimal sufficient sets of chips have size 7, and hence can be distributed to games having $p \ge 4$ players.  \bbox

Given these results, one strategy for survival is for the group to locate one of the minimal sufficient sets in their distribution; if they find neither, its game over. If the set 
\begin{tikzpicture}[scale=0.8, baseline=-1.5mm]
\filldraw[color=blue!60, fill=blue!5, very thick](-1,0) circle (0.3) node {3};
\filldraw[color=red!60, fill=red!5, very thick](-0.2,0) circle (0.3) node {3};
\filldraw[color=green!60, fill=green!5, very thick](0.6,0) circle (0.3) node {1};
\end{tikzpicture}
has been allotted, then they should make the three Rule 1 exchanges 

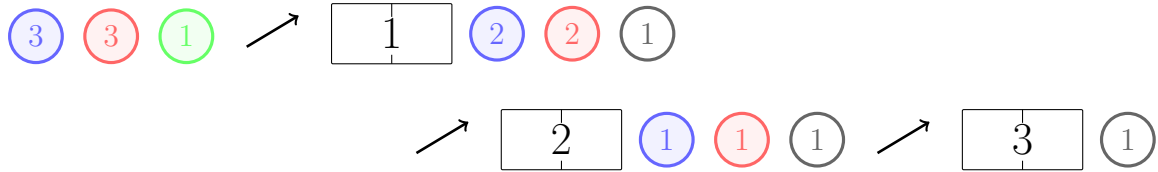
\displaystyled 
\begin{tikzpicture}
	\filldraw[color=blue!60, fill=blue!5, very thick] (-1,0) circle (0.35) node {3};
	\filldraw[color=red!60, fill=red!5, very thick](0,0) circle (0.35) node {3};
	\filldraw[color=green!60, fill=green!5, very thick](1,0) circle (0.35) node {1};
	\end{tikzpicture}
	\begin{tikzpicture}[node distance=2cm]
	\pgfsetlinewidth{1pt}
	\node (A) at (0, 0) {};
	\node (B) at (1, 0.6) {};
	\draw [->] (A) edge (B);
	\end{tikzpicture}
	\begin{tikzpicture}
	\draw[rounded corners=0.5, line hidden] (\doffset,\ndscaleovtwo) rectangle (\domright,\dscaleovtwo);
	\node at (\dommid,0) 
	{\Large{$1$}};
	\draw (\dommid,\dscaleovtwo) -- (\dommid,\legup);
	\draw (\dommid,\ndscaleovtwo) -- (\dommid,\legdown);
	\filldraw[color=blue!60, fill=blue!5, very thick] (-1,0) circle (0.35) node {2};
	\filldraw[color=red!60, fill=red!5, very thick](0,0) circle (0.35) node {2};
	\draw[color=black!60, very thick](1,0) circle (0.35) node {1}; 
	\end{tikzpicture}
\vspace{0.5cm}

\hspace{5cm}
	\begin{tikzpicture}[node distance=2cm]
	\pgfsetlinewidth{1pt}
	\node (A) at (0, 0) {};
	\node (B) at (1, 0.6) {};
	\draw [->] (A) edge (B);
	\end{tikzpicture}
	\begin{tikzpicture}
	\draw[rounded corners=0.5, line hidden] (\doffset,\ndscaleovtwo) rectangle (\domright,\dscaleovtwo);
	\node at (\dommid,0) 
	{\Large{$2$}};
	\draw (\dommid,\dscaleovtwo) -- (\dommid,\legup);
	\draw (\dommid,\ndscaleovtwo) -- (\dommid,\legdown);
	\filldraw[color=blue!60, fill=blue!5, very thick] (-1,0) circle (0.35) node {1};
	\filldraw[color=red!60, fill=red!5, very thick](0,0) circle (0.35) node {1};
	\draw[color=black!60, very thick](1,0) circle (0.35) node {1}; 
	\end{tikzpicture}
	\begin{tikzpicture}[node distance=2cm]
	\pgfsetlinewidth{1pt}
	\node (A) at (0, 0) {};
	\node (B) at (1, 0.6) {};
	\draw [->] (A) edge (B);
	\end{tikzpicture}
		\begin{tikzpicture}
	\draw[rounded corners=0.5, line hidden] (\doffset,\ndscaleovtwo) rectangle (\domright,\dscaleovtwo);
	\node at (\dommid,0) 
	{\Large{$3$}};
	\draw (\dommid,\dscaleovtwo) -- (\dommid,\legup);
	\draw (\dommid,\ndscaleovtwo) -- (\dommid,\legdown);
	\draw[color=black!60, , very thick] (-1,0) circle (0.35) node {1};
	\end{tikzpicture}

\noindent and if allotted 
\begin{tikzpicture}[scale=0.8, baseline=-1.5mm]
\filldraw[color=blue!60, fill=blue!5, very thick](-1,0) circle (0.3) node {3};
\filldraw[color=red!60, fill=red!5, very thick](-0.2,0) circle (0.3) node {2};
\filldraw[color=green!60, fill=green!5, very thick](0.6,0) circle (0.3) node {2};
\end{tikzpicture},
then they should make the three Rule 1 exchanges 

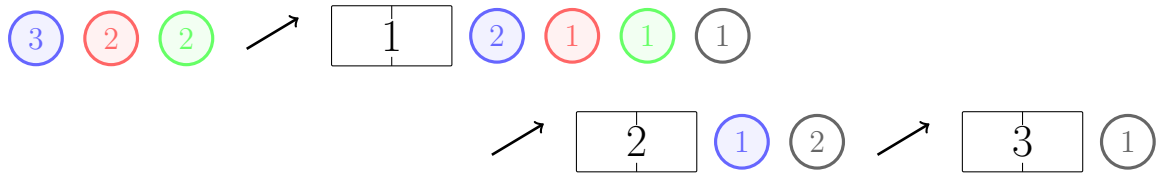
\displaystyled 
\begin{tikzpicture}
\filldraw[color=blue!60, fill=blue!5, very thick] (-1,0) circle (0.35) node {3};
\filldraw[color=red!60, fill=red!5, very thick](0,0) circle (0.35) node {2};
\filldraw[color=green!60, fill=green!5, very thick](1,0) circle (0.35) node {2};
\end{tikzpicture}
\begin{tikzpicture}[node distance=2cm]
\pgfsetlinewidth{1pt}
\node (A) at (0, 0) {};
\node (B) at (1, 0.6) {};
\draw [->] (A) edge (B);
\end{tikzpicture}
\begin{tikzpicture}
\draw[rounded corners=0.5, line hidden] (\doffset,\ndscaleovtwo) rectangle (\domright,\dscaleovtwo);
\node at (\dommid,0) 
{\Large{$1$}};
\draw (\dommid,\dscaleovtwo) -- (\dommid,\legup);
\draw (\dommid,\ndscaleovtwo) -- (\dommid,\legdown);
\filldraw[color=blue!60, fill=blue!5, very thick] (-1,0) circle (0.35) node {2};
\filldraw[color=red!60, fill=red!5, very thick](0,0) circle (0.35) node {1};
\filldraw[color=green!60, fill=green!5, very thick](1,0) circle (0.35) node {1};
\draw[color=black!60, very thick](2,0) circle (0.35) node {1}; 
\end{tikzpicture}
\vspace{0.5cm}

\hspace{6cm}
\begin{tikzpicture}[node distance=2cm]
\pgfsetlinewidth{1pt}
\node (A) at (0, 0) {};
\node (B) at (1, 0.6) {};
\draw [->] (A) edge (B);
\end{tikzpicture}
\begin{tikzpicture}
\draw[rounded corners=0.5, line hidden] (\doffset,\ndscaleovtwo) rectangle (\domright,\dscaleovtwo);
\node at (\dommid,0) 
{\Large{$2$}};
\draw (\dommid,\dscaleovtwo) -- (\dommid,\legup);
\draw (\dommid,\ndscaleovtwo) -- (\dommid,\legdown);
\filldraw[color=blue!60, fill=blue!5, very thick] (-1,0) circle (0.35) node {1};
\draw[color=black!60, very thick](0,0) circle (0.35) node {2}; 
\end{tikzpicture}
\begin{tikzpicture}[node distance=2cm]
\pgfsetlinewidth{1pt}
\node (A) at (0, 0) {};
\node (B) at (1, 0.6) {};
\draw [->] (A) edge (B);
\end{tikzpicture}
\begin{tikzpicture}
\draw[rounded corners=0.5, line hidden] (\doffset,\ndscaleovtwo) rectangle (\domright,\dscaleovtwo);
\node at (\dommid,0) 
{\Large{$3$}};
\draw (\dommid,\dscaleovtwo) -- (\dommid,\legup);
\draw (\dommid,\ndscaleovtwo) -- (\dommid,\legdown);
\draw[color=black!60, , very thick] (-1,0) circle (0.35) node {1};
\end{tikzpicture}

\inlinedom
When either of these exchanges are completed, the group should then fire up the $D^3$
machine, making the moves specified in Theorem \ref{thm:D3}. How much work is in store from them is answered in the next section, where we find that, depending on $I$, some of the use of the $D^3$ machine may be replaced by a more efficient, and thus faster, plan. 

\section{Complexity} \label{sec:complexity}
Here we consider the issue of computing $e(p)$, the number of exchanges required for a group of size $p \ge 4$ to achieve success. As  $e(p)$ depends on the intial distribution of chips, we focus on the two extreme cases, the least favorable one where the group starts from a minimal sufficient set and the remaining chips of a single color, and the most favorable case, where $I$ contains as many `full sets' 
\begin{tikzpicture}[scale=0.8, baseline=-1.5mm]
\filldraw[color=blue!60, fill=blue!5, very thick](-1,0) circle (0.3) node {1};
\filldraw[color=red!60, fill=red!5, very thick](-0.2,0) circle (0.3) node {1};
\filldraw[color=green!60, fill=green!5, very thick](0.6,0) circle (0.3) node {1};
\end{tikzpicture}
as possible. The order of growth of $e(p)$ in $p$ tells us how complex the problem is to resolve in practice. 

For the first case, start with either minimal sufficient set of Theorem \ref{thm:minsuf}, both of size 7, and the remaining $q=2p-7$ chips of a single color. It will not be difficult to see from what follows that the amount of work does not depend on the color of these remaining chips, and we take them to be
\begin{tikzpicture}[scale=0.8, baseline=-1.5mm]
\filldraw[color=blue!60, fill=blue!5, very thick](-1,0) circle (0.3) node {\small{q}};
\end{tikzpicture}.

As was shown at the end of the previous section, 3 exchanges are needed to convert either minimal sufficient set into 
\begin{tikzpicture}[scale=0.8, baseline=-1.5mm]
\draw[rounded corners=0.5, line hidden] (\doffset,\ndscaleovtwo) rectangle (\domright,\dscaleovtwo);
\node at (\dommid,0) 
{$3$};
\end{tikzpicture}
\begin{tikzpicture}[scale=0.8, baseline=-1.5mm]
\draw[color=black!60, thick](-1,0) circle (0.3) node {1};
\end{tikzpicture}.
Next, by Theorem \ref{thm:D3} for $r=1$, four exchanges are needed to obtain
\begin{tikzpicture}[scale=0.8, baseline=-1.5mm]
\draw[rounded corners=0.5, line hidden] (\doffset,\ndscaleovtwo) rectangle (\domright,\dscaleovtwo);
\node at (\dommid,0) 
{$3$};
\end{tikzpicture}
\begin{tikzpicture}[scale=0.8, baseline=-1.5mm]
\filldraw[color=blue!60, fill=blue!5,thick](-1,0) circle (0.3) node {\small{q}};
\draw[color=black!60, thick](-0.2,0) circle (0.3) node {2};
\end{tikzpicture}
and a single further Rule 1 exchange to achieve
\begin{tikzpicture}[scale=0.8, baseline=-1.5mm]
\draw[rounded corners=0.5, line hidden] (\doffset,\ndscaleovtwo) rectangle (\domright,\dscaleovtwo);
\node at (\dommid,0) 
{$4$};
\end{tikzpicture}
\begin{tikzpicture}[scale=0.8, baseline=-1.5mm]
\filldraw[color=blue!60, fill=blue!5,thick](-1,0) circle (0.3) node {\tiny{q-1}};
\draw[color=black!60, thick](-0.2,0) circle (0.3) node {1};
\end{tikzpicture},
for a total of five. As these last five exchanges need to be run $p-3$ times to achieve
\begin{tikzpicture}[scale=0.8, baseline=-1.5mm]
\draw[rounded corners=0.5, line hidden] (\doffset,\ndscaleovtwo) rectangle (\domright,\dscaleovtwo);
\node at (\dommid,0) 
{$p$};
\end{tikzpicture},
\setcounter{equation}{\eqhinge}
the total number of exchanges required is
\begin{align} \label{eq:e(p).mult3}
e(p)=5(p-3)+3=5p-12 \qmq{for $p \ge 4$.}
\end{align}

In particular, the complexity of the problem is linear in its size. Even so, this plan may be cumbersome to implement. 
For instance, the group of \p\ players considered in the introduction who were given \minutes\ minutes to play, having only the worst case set $I$ at their disposal, are required to make $e(\p)=\minabc$ exchanges. Even if exchanges could be made every \exrate\ seconds with no errors or delays, completing this many would take 19 \signore{\marginpar{\minsminabc\ }} minutes. Exchanges at this pace may be even more unrealistic in the version of the game without a joker, where the players need to keep track of the current least frequent color.

On the other hand, the group may be given an initial distribution of chips with many copies of the `full' set 
\begin{tikzpicture}[scale=0.8, baseline=-1.5mm]
\filldraw[color=blue!60, fill=blue!5, very thick](-1,0) circle (0.3) node {1};
\filldraw[color=red!60, fill=red!5, very thick](-0.2,0) circle (0.3) node {1};
\filldraw[color=green!60, fill=green!5, very thick](0.6,0) circle (0.3) node {1};
\end{tikzpicture}
that can be traded for a domino and a joker in one economical move. To find the maximum advantage that would be gained by making use of full sets, consider the case where $p$ is some multiple of three, so that we may write $2p=3m$,  and consider the `best case' initial distribution that has the maximum number of full sets, 
\begin{tikzpicture}[scale=0.8, baseline=-1.5mm]
\filldraw[color=blue!60, fill=blue!5, very thick](-1,0) circle (0.3) node {$m$};
\filldraw[color=red!60, fill=red!5, very thick](-0.2,0) circle (0.3) node {$m$};
\filldraw[color=green!60, fill=green!5, very thick](0.6,0) circle (0.3) node {$m$};
\end{tikzpicture};
other cases where the maximum number of `full' sets are initially distributed is left for the reader. Here, for non-triviality, we must take $p \ge 6$, the smallest multiple of 3 that is at least 4. 

Applying $m$ Rule 1 exchanges, followed by Theorem \ref{thm:Xtox}, and the observation that $m$ must be even, yields
\displaystyled
\newcommand{\circradtwo}{0.45}
\newcommand{\circspacetwo}{0.95}
\FPeval{\cicspacetwot}{\circspacetwo-1}
\FPeval{\cicspacethreet}{2*\circspacetwo-1}
\begin{center}
	\begin{tikzpicture}
	\filldraw[color=blue!60, fill=blue!5, very thick] (-1,0) circle (\circradtwo) node {\Large{$m$}};
	\filldraw[color=red!60, fill=red!5, very thick](\cicspacetwot,0) circle (\circradtwo) node {\Large{$m$}};
	\filldraw[color=green!60, fill=green!5, very thick](\cicspacethreet,0) circle (\circradtwo) node {\Large{$m$}};
	\end{tikzpicture}
	\begin{tikzpicture}[node distance=2cm]
	\pgfsetlinewidth{1pt}
	\node (A) at (0, 0) {};
	\node (B) at (1, 0.6) {};
	\draw [->] (A) edge (B);
	\end{tikzpicture}
	\begin{tikzpicture}
	\draw[rounded corners=0.5, line hidden] (\doffset,\ndscaleovtwo) rectangle (\domright,\dscaleovtwo);
	\node at (\dommid,0) 
	{\Large{$m$}};
	\draw (\dommid,\dscaleovtwo) -- (\dommid,\legup);
	\draw (\dommid,\ndscaleovtwo) -- (\dommid,\legdown);
	\draw[color=black!60, very thick](-1,0) circle (\circradtwo) node {\Large{$m$}};
	\end{tikzpicture}
	\begin{tikzpicture}[node distance=2cm]
	\pgfsetlinewidth{1pt}
	\node (A) at (0, 0) {};
	\node (B) at (1, 0.6) {};
	\draw [->] (A) edge (B);
	\end{tikzpicture}
	\begin{tikzpicture}
	\draw[rounded corners=0.5, line hidden] (\doffset,\ndscaleovtwo) rectangle (\domright,\dscaleovtwo);
	\node at (\dommid,0) 
	{\small{$m+d(m)$}};
	\draw (\dommid,\dscaleovtwo) -- (\dommid,\legup);
	\draw (\dommid,\ndscaleovtwo) -- (\dommid,\legdown);
	\draw[color=black!60, very thick](-1,0) circle (\circradtwo) node {\small{$\phi(m)$}};
	\node[] at (0,0) {$=$};
	\end{tikzpicture}
    \begin{tikzpicture}
		\draw[rounded corners=0.5, line hidden] (\doffset,\ndscaleovtwo) rectangle (\domright,\dscaleovtwo);
	\node at (\dommid,0) 
	{\small{$p-1$}};
	\draw (\dommid,\dscaleovtwo) -- (\dommid,\legup);
	\draw (\dommid,\ndscaleovtwo) -- (\dommid,\legdown);
		\draw[color=black!60, very thick](-1,0) circle (\circradtwo) node {2};
			\node[] at (0.5,0) {(\eqhingpltwo)};
	\end{tikzpicture}
\end{center}

By Theorem \ref{thm:Xtox} the number of exchanges required to reach this state is $m+d(m)=p-1$.  
To continue we must invoke Rule 2, which is possible as $p-1 \ge 3$. Then, to obtain the final domino we apply the exchanges outlined in Theorem \ref{thm:Xtox}, yielding
\begin{center}
	\begin{tikzpicture}
	\draw[rounded corners=0.5, line hidden] (\doffset,\ndscaleovtwo) rectangle (\domright,\dscaleovtwo);
	\node at (\dommid,0) 
	{$p-1$};
	\draw (\dommid,\dscaleovtwo) -- (\dommid,\legup);
	\draw (\dommid,\ndscaleovtwo) -- (\dommid,\legdown);
	\draw[color=black!60, very thick](-1,0) circle (\circradtwo) node {2};
	\end{tikzpicture}
	\begin{tikzpicture}[node distance=2cm]
	\pgfsetlinewidth{1pt}
	\node (A) at (0, 0) {};
	\node (B) at (1, 0.6) {};
	\draw [->] (A) edge (B);
	\end{tikzpicture}
	\begin{tikzpicture}
	\draw[rounded corners=0.5, line hidden] (\doffset,\ndscaleovtwo) rectangle (\domright,\dscaleovtwo);
	\node at (\dommid,0) 
	{$p-4$};
	\draw (\dommid,\dscaleovtwo) -- (\dommid,\legup);
	\draw (\dommid,\ndscaleovtwo) -- (\dommid,\legdown);
	\draw[color=black!60, very thick](-1,0) circle (\circradtwo) node {9};
	\end{tikzpicture}
	\begin{tikzpicture}[node distance=2cm]
	\pgfsetlinewidth{1pt}
	\node (A) at (0, 0) {};
	\node (B) at (1, 0.6) {};
	\draw [->] (A) edge (B);
	\end{tikzpicture}
	\begin{tikzpicture}
	\draw[rounded corners=0.5, line hidden] (\doffset,\ndscaleovtwo) rectangle (\domright,\dscaleovtwo);
	\node at (\dommid,0) 
	{\small{$p$-4+$d$(9)}};
	\draw (\dommid,\dscaleovtwo) -- (\dommid,\legup);
	\draw (\dommid,\ndscaleovtwo) -- (\dommid,\legdown);
	\draw[color=black!60, very thick](-1,0) circle (\circradtwo) node {\small{$\phi(9)$}};
	\node[] at (0,0) {$=$};
	\end{tikzpicture}
	\begin{tikzpicture}
	\draw[rounded corners=0.5, line hidden] (\doffset,\ndscaleovtwo) rectangle (\domright,\dscaleovtwo);
	\node at (\dommid,0) 
	{\small{$p$}};
	\draw (\dommid,\dscaleovtwo) -- (\dommid,\legup);
	\draw (\dommid,\ndscaleovtwo) -- (\dommid,\legdown);
	\draw[color=black!60, very thick](-1,0) circle (\circradtwo) node {1};
	\node[] at (0.5,0) {(\eqhingplthree)};
	\end{tikzpicture}
\end{center}

The  Rule 2 exchange is followed by $d(9)=4$ applications of Rule 1, by Theorem \ref{thm:Xtox}, for a total of 5 exchanges. Hence, adding in the cost of $p-1$ for (\eqhingpltwo) we obtain
\begin{align*}
e(p)= p+4.
\end{align*}
The order $e(p)$ is still linear in $p$, but the multiplier has dropped from 5 to 1, which is quite a savings. Returning to our  group of $p=$\p\ players, satisfying $2p=3m$ for $m=32$, we see that with this most favorable $I$, the work decreases from \minabc\ exchanges to \abc, and time required drops from 19 minutes to just over 4,
at \exrate\ seconds per exchange.\signore{\marginpar{\minsminabc\,\minsabc}} For a game fixed at \minutes\ minutes, the facilitator can make the game anywhere between near impossible to completely easy.

\inlinedom A rough computation shows one can expect this type of savings quite generally in the case where $I$ is most favorable, that is, when $2p=3m+r, r \in \{0,1,2\}$ and $m$ full sets are distributed.  To avoid the handling of cumbersome `boundary cases' we let $m \ge 3$. Making the first two transformations in (\eqhingpltwo) yields $m+d(m)$ dominoes, costing that many exchanges, and $\phi(m) \ge 1$ jokers. 
Using only the subset
\begin{tikzpicture}[scale=0.8, baseline=-1.5mm]
\draw[rounded corners=0.5, line hidden] (\doffset,\ndscaleovtwo) rectangle (\domright,\dscaleovtwo);
\node at (\dommid,0) 
{$3$};
\end{tikzpicture}
\begin{tikzpicture}[scale=0.8, baseline=-1.5mm]
\draw[color=black!60, thick](-1,0) circle (0.3) node {1};
\end{tikzpicture},
the $D^3$ machine can be used to produce the remaining needed $p-(m+d(m))$ dominoes, by Theorem \ref{thm:D3}, at 9 exchanges per use. Hence we obtain the upper bound
\begin{multline*}
e(p) \le  m+d(m)+9(p-m-d(m))=9p-8(m+d(m)) \le 9p-8\left(3m/2-1\right)\\
= 9p-8\left(p-r/2-1\right)=p+4r+8 \le p+16, 
\end{multline*}
using $d(m) \ge m/2-1$ for the first inequality, that $m=(2p-r)/3$, and then $r \le 2$
for the final inequality.  Comparing this upper bound with the precise value of $e(p)$ given in \eqref{eq:e(p).mult3} in one special case gives some evidence that this `rough' bound isn't very off from the truth in general.  

\section{Conclusion}

The game of chips, dominoes and survival is an exciting team building exercise that enforces that the group learn to cooperate at a high level. The existence of two small sets sufficient for stoking the $D^3$ machine and, guaranteeing the survival of a group of any size, given sufficient time to make exchanges, may be somewhat surprising. The group can take advantage of the full sets in their initial distribution, but only to the extent that a facilitator can control. In cases where the least favorable set of chips is distributed, the abundance of chips in larger games may contribute more to disorder than to quick resolution. 

This game and its analysis, though simple on the one hand, contains glimmerings of a much wider landscape, including cooperative game theory, logical reasoning, number theory and the computation of algorithmic complexities. 
We hope readers not only have fun playing and facilitating this game, but also in dreaming up extensions with higher level pieces and more complex exchange rules, and seeing what new and surprising properties these survival games may have in store.

\end{document}